\documentclass{article}

\usepackage{amsmath}
\usepackage{url}

\usepackage{amssymb}
\usepackage{amsmath}
\usepackage{amsthm}
\usepackage{amsfonts}
\usepackage{mathrsfs}
\usepackage{makeidx}
\usepackage{multicol}

\def \qed{\hfill $\Box$}

\begin{document}

\title{On Picard Type Theorems and Entire Solutions of Differential Equations}
\author{Bao Qin Li and  Liu Yang }

\date{}
\maketitle

\bigskip
\noindent {\it Abstract}. We give a connection between the Picard
type theorem of Polya-Saxer-Milliox and characterization of entire
solutions of a differential equation and then their higher
dimensional extensions, which leads further results on both
(ordinary and partial) differential equations and Picard type
theorems.



\vskip.5in
\bigskip

In the recent paper \cite{Li1}, a connction/equivalence was
established between Picard's theorem in complex analysis and
characterization of entire solutions of an ordinary differential
equation, which was further extended to partial differential
equations (see \cite{Li1} for the details). Recall that Picard's
theorem asserts that an entire function, i.e., a complex-valued
function differentiable in the complex plane $\textbf{C}$, omitting
two complex numbers must be constant. It also implies, by a simple
transform, the meromorphic version of the theorem that a meromorphic
function in $\textbf{C}$ omitting three distinct values must be
constant. Picard's theorem is among the most striking results in
complex analysis and plays a decisive role in the development of the
theory of entire and meromorphic functions and other applications.
We refer to \cite{Se} for an exposition (the history, methods, and
references) of the theorem and \cite{ah}, \cite{bo}, \cite{da},
\cite{fu}, \cite{Le}, \cite{zh}, etc. for various proofs.

We call a theorem a Picard type theorem if it asserts that an entire
function under certain conditions must reduce to constant. The
observation in \cite{Li1} suggests a heuristic: A Picard type
theorem is apt to imply a characterization of entire solutions of a
differential equation, and vice versa. (This may also be extended to
the meromorphic case.)

Inspired by this heuristic, we present, in this article, a
connection between the well-known Picard type theorem of
Polya-Saxer-Milliox (see \cite{Sa} and \cite{Mi}; see also
\cite{Ha}) and a characterization of entire solutions of an ordinary
differential equation, which provides a different approach to the
existing results and leads further results on both (ordinary and
partial) differential equations and Picard type theorems (cf.
below).

\bigskip
\noindent {\bf Theorem A} (Polya-Saxer-Millioux) {\it  An entire
function $f$ must be constant if $f$ omits $0$ and $f'$ omits $1$ in
$\textbf{C}$.}

\medskip
We note that Theorem A has undergone various extensions over the
years; especially it also holds for meromorphic functions (see
\cite{Ha}). In this paper we however only consider entire
functions, which serves our main purpose for the above mentioned
connections, although Theorem 1 below holds also for
meromorphic functions, which is clear from the proof given below.

We also note, by a simple linear transform, that the values $0$ and
$1$ in Theorem A can be replaced by any two complex numbers $a$ and
$b\not=0$, respectively; but, the assumption $b\not=0$ cannot be
dropped, as seen from the nonconstant entire function $f=e^z$, for
which both $f$ and $f'$ omit $0$.

We first give the characterization of entire solutions of the differential equation in
the following

\bigskip
\noindent {\bf Theorem 1.} {\it Let $a(z)$ be an entire function in
$\textbf{C}$ and let $P(z_1, z_2)$ be an entire function in
$\textbf{C}^2$ divisible by $(z_1-c)(z_2-d)$ in the ring of entire functions, where $c, d$ with $d\not=0$ are complex numbers. Then an
entire solution $f$ in $\textbf{C}$ of the differential equation
\begin{equation}\label{eq1}
f'+a(z)P(f, f')=0
\end{equation}
must be constant.
}

\medskip
We will see that Theorem 1 implies Theorem A immediately; as a matter of fact, Theorem 1 and Theorem A are equivalent in the sense that one implies the other.

\bigskip
\noindent {\bf Theorem 1 $\Longrightarrow$ Theorem A}. Since $f$ satisfies that $f\not=0$ and $f'\not=1$, the function
$a(z):={f'\over f(f'-1)}$ is entire. Clearly, $f'-a(z)f(f'-1)=0$,
which is an equation of the form (\ref{eq1}) with
$P(z_1, z_2)=z_1(z_2-1)$. Thus, $f$ must be constant by Theorem 1.
\qed

\bigskip
\noindent {\bf Theorem A $\Longrightarrow$ Theorem 1}. Since $P$ is
divisible by $(z_1-c)(z_2-d)$ in the ring of entire functions, we
can write $P(z_1, z_2)=(z_1-c)(z_2-d)g(z_1, z_2)$, where $g$ is an
entire function in $\textbf{C}^2$. It follows from (\ref{eq1}) that
\begin{equation}\label{zero}
f'(z)=-a(z)(f(z)-c)(f'(z)-d)g(f(z), f'(z)).
\end{equation}
It is clear that $f'$ cannot assume $d$, since otherwise, the right
hand side of (\ref{zero}) is $0$ while the left hand side is
nonzero, a contradiction. It is also easy to see that $f$ cannot
assume $c$, since otherwise the right hand side of (\ref{zero})
would have a zero (coming from a zero of $f-c$) with multiplicity
strictly greater than that of the same zero of the left hand side
(due to the derivative, which decreases the multiplicity), which is
absurd. Thus, $f$ omits $c$ and $f'$ omits $d$. By Theorem A,
$F={f-c\over d}$ and thus $f$ must be constant. \qed

\medskip
Study of solutions of differential equations in complex variables
has a long history. The connection between Picard type theorems and
characterizations of entire solutions of differential equations may
help discover further results in both directions, as seen below.

Theorem 1, equivalent to the Picard type theorem of
Polya-Saxer-Millioux but in the form on solutions of the
differential equation, leads naturally consideration of the related higher
order differential equations
\begin{equation}\label{ode1}
f^{(n)}+a(z)P(f, f'))=0
\end{equation}
for $n\ge 2$, in view of the known fact that Theorem A still holds
with $f^{(n)}$ replacing $f'$ (see e.g. \cite{Ha}). It is tempting
to expect a result similar to Theorem 1 would hold for the equation
(\ref{ode1}). However, entire solutions of (\ref{ode1}) are not
necessarily constant any more. In fact, when $n\ge 2$ the equation
$f^{(n)}+a(z)f(f'-1)=0$ has a nonconstant entire solutions $f(z)=z$. We will see that there is a generality
behind this and a complete characterization of entire
solutions can be given, which allows extensions to more general partial
differential equations and can yield further results
on Picard type theorems. We consider

\begin{equation}\label{pde}
\sum\limits_{|\alpha|=1}^ma_{\alpha}\frac{\partial^{|\alpha|}f_{z_j}(z)}{\partial^{\alpha_1}z_1\cdots\partial^{\alpha_n}z_n}
+a(z)P(f(z), f_{z_j}(z))=0
\end{equation}
where $z=(z_1, z_2, \cdots, z_n)$ in $\textbf{C}^n$,
$\alpha=(\alpha_1,\cdots, \alpha_n)$ is a multi-index with
$\alpha_j\ge 0$ and  $|\alpha|=\alpha_1+\cdots+\alpha_n$,
$a_{\alpha}$'s are polynomials in $\textbf{C}^n$ and
$f_{z_j}={\partial f\over \partial z_j}$ ($1\le j\le n$).

\bigskip
\noindent {\bf Theorem 2.} {\it Let $a(z)$ be a nonzero entire
function in $\textbf{C}^n$ and let $P(x, y)$ be a nonzero entire
function in $\textbf{C}^2$ divisible by $(x-c)(y-d)$ in the ring of
entire functions in $\textbf{C}^2$ for two complex numbers $c, d$
with $d\not=0$. Then an entire solution $f\not\equiv c$ in
$\textbf{C}^n$ of the partial differential equation (\ref{pde}) is
given by $f(z)=dz_j+\phi$, where $\phi$ is an entire function in
$z_1,\cdots, z_{j-1}, z_{j+1},\cdots,z_n$. In particular,
$f(z)=dz+A$ if $n=1$, where $A$ is a constant.}

\medskip
In proving Theorem 2, we will utilize the following known properties
of
$$m(r, f):=\int_{S_n(r)}\log^+|f|\eta_n$$
for a nonconstant entire function $f$ in $\textbf{C}^n$, where
$\log^+ x=\max\{0, \log x\}$ and $\eta_n$ is the usual positive
volume form on the sphere $S_n(r):=\{z\in\textbf{C}^n: |z|=r\}$
normalized so that the total volume of the sphere is $1$ (see e.g.
\cite{Vi}; cf. also \cite{Li1}):

  \medskip
(I) $m(r, 1/f)\le m(r, f)+O(1)$, which follows from Jensen's formula
or the Nevanlinna first fundamental theorem in several complex
variables;

(II) $m(r,
\frac{\frac{\partial^{|\alpha|}f(z)}{\partial^{\alpha_1}z_1\cdots\partial^{\alpha_n}z_n}}{f})
=S(r, f)$, where $S(r, f)$ denotes a quantity satisfying that $S(r,
f)\doteq O\{\log (rm(r, f))\}$ as $r\to \infty$, which follows from
the logarithmic derivative lemma in several complex variables,
where the symbol $\doteq$ means that the equality holds outside a
set of $r$ of finite Lebesgue measure.

\bigskip
\noindent {\bf Proof of Theorem 2}. If
$f_{z_j}\equiv d$, the theorem already holds. We assume in the following that $f_{z_j}\not\equiv d$
and will derive a contradiction.

By the assumption on $P$, we can write $P(x, y)=(x-c)(y-d)g(z),$
where $g$ is a nonzero entire function in $\textbf{C}^2$. Thus, the
given equation can be written as
\begin{eqnarray}\label{pde-1}
& &\sum\limits_{|\alpha|=1}^ma_{\alpha}\frac{\partial^{|\alpha|}f_{z_j}(z)}{\partial^{\alpha_1}z_1\cdots\partial^{\alpha_n}z_n}
\nonumber\\
& &=b(z)(f(z)-c)(f_{z_j}(z)-d),
\end{eqnarray}
where $b(z)=-a(z)g(f(z), f_{z_j}(z))$. We now write (\ref{pde-1}) as
\begin{eqnarray*}
\frac{\sum\limits_{|\alpha|=1}^ma_{\alpha}\frac{\partial^{|\alpha|}f_{z_j}(z)}{\partial^{\alpha_1}z_1
\cdots\partial^{\alpha_n}z_n}}{f_{z_j}(z)-d}=b(z)(f(z)-c).
\end{eqnarray*}
and
\begin{eqnarray*}
\frac{\sum\limits_{|\alpha|=1}^ma_{\alpha}\frac{\partial^{|\alpha|}f_{z_j}(z)}{\partial^{\alpha_1}z_1
\cdots\partial^{\alpha_n}z_n}}{f(z)-c}=b(z)(f_{z_j}(z)-d)
\end{eqnarray*}
in view of the fact that $f\not\equiv c$.  We then obtain by
Property (II) that
$$m(r, b(f-c))=S(r, f_{z_j})=S(r, f)$$
and
$$ m(r, b(f_{z_j}-d)=S(r, f),$$
which implies, by Property (I), that
\begin{eqnarray*}
& &m(r, {f_{z_j}-d\over f-c})=m(r, {b(f_{z_j}-d))\over b(f-c)})\\
& &\le m(r, b(f_{z_j}-d))+m(r, {1\over b(f-c)})\\
& &\le m(r, b(f_{z_j}-d))+m(r, b(f-c))+O(1)=S(r, f).
\end{eqnarray*}
We then deduce that
\begin{eqnarray*}
& &m(r, {1\over f-c})\le m(r, {d\over f-c})+O(1)\\
& &\le m(r, {d-f_{z_j}\over f-c})+m(r, {f_{z_j}\over f-c})+O(1)=S(r, f).
\end{eqnarray*}
Changing the equation (\ref{pde-1}) to
\begin{eqnarray*}
\frac{\sum\limits_{|\alpha|=1}^ma_{\alpha}\frac{\partial^{|\alpha|}f_{z_j}(z)}{\partial^{\alpha_1}z_1
\cdots\partial^{\alpha_n}z_n}}{f_{z_j}(z)-d}{1\over f-c}=b(z)
\end{eqnarray*}
we obtain that $m(r, b)=S(r, f)$. We then change the equation
(\ref{pde-1}) to
\begin{eqnarray*}
\frac{\sum\limits_{|\alpha|=1}^ma_{\alpha}\frac{\partial^{|\alpha|}f_{z_j}(z)}{\partial^{\alpha_1}z_1
\cdots\partial^{\alpha_n}z_n}}{f_{z_j}(z)-d}{1\over b}=f-c,
\end{eqnarray*}
from which and Properties (I) and (II) again we deduce that
\begin{eqnarray*}
& &m(r, f)=m(r, f-c)+O(1) \\
& &\le m(r, \frac{\sum\limits_{|\alpha|=1}^ma_{\alpha}\frac{\partial^{|\alpha|}f_{z_j}(z)}{\partial^{\alpha_1}z_1
\cdots\partial^{\alpha_n}z_n}}{f_{z_j}(z)-d})+m(r, {1\over b})\\
& &=S(r, f)\doteq O\{\log
(rm(r, f))\}
\end{eqnarray*}
and then that $m(r, f)\doteq O(\log r)$, which implies that
$f$ is a polynomial. Then, the left hand side of (\ref{pde-1}), as a
linear combination of polynomials, must be a polynomial, which
implies that the entire function $b(z)$, as a quotient of two
polynomials, must be a polynomial, too. But, it is evident that the
degree of the left hand side of (\ref{pde-1}) (due to higher order
derivatives) is strictly lower than that of the left hand side,
which is impossible. This completes the proof. \qed

\medskip
The above proof of Theorem 2 can be pushed over to even more general
situation where $f_{z_j}$ is replaced by a partial derivative
$\frac{\partial^{|I|}f(z)}{\partial^{i_1}z_1\cdots\partial^{i_n}z_n}$
of any order or a sum of the form
\begin{equation}\label{Df}
Df=\sum_{|I|=1}^ma_I\frac{\partial^{|I|}f(z)}{\partial^{i_1}z_1\cdots\partial^{i_n}z_n}
\end{equation}\label{pde2}
with $I=(i_1, \cdots, i_n)$ and $a_I$s being polynomials in
$\textbf{C}^n$ (or any expression formed by $f$ and finitely many
partial derivatives of $f$ of any orders so that the proof can go
through). We include the following theorem for partial differential
equations
\begin{equation}\label{pde2}
\sum\limits_{|\alpha|=1}^ma_{\alpha}\frac{\partial^{|\alpha|}Df(z)}{\partial^{\alpha_1}z_1\cdots\partial^{\alpha_n}z_n}
+a(z)P(f(z), Df(z))=0.
\end{equation}

\medskip \noindent {\bf Theorem 3.} {\it  Let $a(z)$ be a nonzero entire function in $\textbf{C}^n$ and let $P(x, y)$ be a nonzero
entire function in $\textbf{C}^2$ divisible by $(x-c)(y-d)$ in the
ring of entire functions in $\textbf{C}^2$ for two complex numbers
$c, d$ with $d\not=0$.  Then an entire function $f\not\equiv c$ in
$\textbf{C}^n$ is a solution of the partial differential equation
(\ref{pde2}) if and only if $f$ is a solution of $Df(z)=d$.}

\medskip
The sufficiency of Theorem 3 is clear in view of (\ref{pde-1}) with
$f_{z_j}$ replaced by $Df$. The proof of the necessity of Theorem 3
is identical to that of Theorem 2 by replacing $f_{z_j}$ there with
$Df$. We thus omit the details.

\medskip
Theorem 3 can yield some old and new Picard type theorems. In particular,
Theorem 3 generalizes the result of Polya-Saxer-Milliox in the
case $n=1$ and thus provides a different approach to the Picard type
theorem.

\medskip \noindent {\bf Corollary 4.} {\it  Suppose that $f$ is an entire function in $\textbf{C}^n$.
If $f$ omits $0$ and $Df$ omits $1$, then $Df$ is identically zero.
}

\medskip \noindent {\bf Examples.} (i) Take $n=1$ and $Df=f^{(m)}$ in Corollary 4, then
the conclusion $Df\equiv 0$ implies, by integration, that $f$ is a
polynomial and thus $f$ must be constant since $f$ omits $0$. Thus,
Corollary 4 gives the Picard type theorem of Polya-Saxer-Milliox
that an entire function $f$ must be constant if $f$ omits $0$ and
$f^{(m)}$ omits $1$.

(ii) Take $n=1$ and $Df=f^{(m)}-f^{(m-1)}$, where $m$ is any
positive integer. Then under the condition of Corollary 4, we must
have $Df=0$, i.e., $f^{(m+1)}-f^{(m)}=0$. Solving this linear
differential equation directly, we obtain that
$f=c_0+c_1z+\cdots+c_{m-1}z^{m-1}+c_me^z$, where $c_j'$s are
constants. But, $f$ does not assume $0$; thus we must have that
$f=c_0$, a constant.

(iii) While $Df$ is identically zero in the conclusion of Corollary
4, the function $f$ itself is not necessarily constant. Take
$f=e^{z_1-z_2}$, a nonconstant entire function in $\textbf{C}^2$.
But, $f$ omits $0$ and $Df:=f_{z_1}+f_{z_2}\equiv 0$ omits $1$,
satisfying the conditions of Corollary 4. The entire functions $f$
are however characterized by the conclusion of Corollary 4, i.e.,
$Df=f_{z_1}+f_{z_2}=0$. In fact, solving this partial differential
equation directly by using its characteristic equations, one obtains
that $f=g(z_1-z_2)$, where $g$ is an entire function in
$\textbf{C}.$ But $f$ does not assume $0$; thus $g=e^h$ for an
entire function $h$. Hence $f=e^{h(z_1-z_2)}.$

\medskip \noindent
{\bf Proof of Corollary 4}. For any fixed $1\le j\le n$, the
function $a(z):=\frac{\frac{\partial Df}{\partial z_j}}{f(Df-1)}$ is
entire. We have that
$$\frac{\partial Df}{\partial
z_j}=af(Df-1).$$  If $a\not\equiv 0$, then by the necessary
condition of Theorem 3, we have that $Df=1$, a contradiction to the
assumption that $Df$ omits $1$. Thus, $a\equiv 0$, i.e., $\frac
{\partial Df}{\partial z_j}=0$. We have this equality for all $1\le
j\le n$, which implies that $Df=C$ is constant. Since $f$ omits $0$
we have that $f=e^g$ for an entire function $g$ in $\textbf{C}^n$.
By the definition of $Df$, it is clear that the equality $Df=C$
becomes $e^gQ=C$, where $Q$ is a polynomial in some partial
derivatives of $g$ with polynomial coefficients. We claim that $Q$
must be identically zero. If not, $g$ then cannot be constant by the
definition of $Q$ and we can write the above equality to
$e^g={C\over Q}$, which implies, by applying Properties I and  II,
that $m(r, e^g)=m(r, {C\over Q})\doteq O(m(r, g))+O(\log r)$ as
$r\to\infty$. But, whenever $g$ is not constant it always holds that
$m(r, g)=o(T(r, e^g))$ by the several complex variable version of
the Clunie's lemma (see \cite{CLY}, p.88). We then have that $m(r,
e^g)=o(m(r, e^g))+O(\log r)$. This is impossible unless $g$ is
constant, a contradiction. This shows the claim, i.e., $Q=0$ and
thus $C=0$. Therefore, $Df$ is identically zero. \qed

\medskip

\bigskip
\noindent\textit{Department of Mathematics and Statistics, Florida
International University, Miami, FL 33199 USA \\ 
libaoqin@fiu.edu}

\bigskip
\noindent\textit{School of Mathematics \& Physics Science and Engineering \\ Anhui University of
Technology,  Ma'anshan, 243032, P.R. China      \\  yangliu20062006@126.com}

\end{document}